       \documentclass[12pt]{amsart}
\usepackage{amsmath}
\usepackage{amsxtra}
\usepackage{amscd}
\usepackage{amsthm}
\usepackage{amsfonts}
\usepackage{amssymb}
\usepackage{eucal}
\usepackage[all]{xypic}

\usepackage{color}

\newtheorem{lem}[subsection]{Lemma}

\theoremstyle{definition}

\newtheorem{defn}[subsection]{Definition}

\newtheorem{remark}[subsection]{Remark}

\newcommand{\thmref}[1]{Theorem~\ref{#1}}

\newcommand{\lemref}[1]{Lemma~\ref{#1}}

\newcommand{\propref}[1]{Proposition~\ref{#1}}
\newcommand{\corref}[1]{Corollary~\ref{#1}}

\newcommand{\exref}[1]{Example~\ref{#1}}

\newcommand{\nc}{\newcommand}
\nc{\renc}{\renewcommand}
\nc{\ssec}{\subsection}
\nc{\sssec}{\subsubsection}
\nc{\on}{\operatorname}

\nc\ol{\overline}
\nc\wt{\widetilde}
\nc\wh{\widehat}
\nc\tboxtimes{\wt{\boxtimes}}

\renc{\d}{{\delta}}
\nc{\Aa}{{\mathbb{A}}}
 \nc{\Gg}{{\mathbb{G}}}  
\nc{\Hh}{{\mathbb{H}}}
 \nc{\Nn}{{\mathbb{N}}}
\nc{\Pp}{{\mathbb{P}}}
\nc{\Rr}{{\mathbb{R}}}
\nc{\BV}{{\mathbb{V}}}
\nc{\BW}{{\mathbb{W}}}
\nc{\Zz}{{\mathbb{Z}}}
\nc{\Qq}{{\mathbb{Q}}}
\nc{\Ss}{{\mathbb{S}}}
\nc{\Cc}{{\mathbb{C}}}
\nc{\Ff}{{\mathbb{F}}}

\nc{\CA}{{\mathcal{A}}}
\nc{\CB}{{\mathcal{B}}}

\nc{\CE}{{\mathcal{E}}}
\nc{\CF}{{\mathcal{F}}}
\nc{\CG}{{\mathcal{G}}}
\nc{\CL}{{\mathcal{L}}}
\nc{\CC}{{\mathcal{C}}}
\nc{\CM}{{\mathcal{M}}}
\def\Mm{\CM}
\nc{\CN}{{\mathcal{N}}}
\nc{\Oo}{{\mathcal{O}}}
\nc{\CP}{{\mathcal{P}}}
\nc{\CQ}{{\mathcal{Q}}}
\nc{\CR}{{\mathcal{R}}}
\nc{\CS}{{\mathcal{S}}}
\nc{\CT}{{\mathcal{T}}}
\nc{\CU}{{\mathcal{U}}}
\nc{\CV}{{\mathcal{V}}}
\nc{\CK}{{\mathcal{K}}}
\nc{\CW}{{\mathcal{W}}}
\nc{\CZ}{{\mathcal{Z}}}

\nc{\cM}{{\check{\mathcal M}}{}}
\nc{\csM}{{\check{\mathcal A}}{}}
\nc{\oM}{{\overset{\circ}{\mathcal M}}{}}
\nc{\obM}{{\overset{\circ}{\mathbf M}}{}}
\nc{\oCA}{{\overset{\circ}{\mathcal A}}{}}
\nc{\obA}{{\overset{\circ}{\mathbf A}}{}}
\nc{\ooM}{{\overset{\circ}{M}}{}}
\nc{\osM}{{\overset{\circ}{\mathsf M}}{}}
\nc{\vM}{{\overset{\bullet}{\mathcal M}}{}}
\nc{\nM}{{\underset{\bullet}{\mathcal M}}{}}
\nc{\oD}{{\overset{\circ}{\mathcal D}}{}}
\nc{\obD}{{\overset{\circ}{\mathbf D}}{}}
\nc{\oA}{{\overset{\circ}{\mathbb A}}{}}
\nc{\op}{{\overset{\bullet}{\mathbf p}}{}}
\nc{\cp}{{\overset{\circ}{\mathbf p}}{}}
\nc{\oU}{{\overset{\bullet}{\mathcal U}}{}}
\nc{\oZ}{{\overset{\circ}{\mathcal Z}}{}}
\nc{\ofZ}{{\overset{\circ}{\mathfrak Z}}{}}
\nc{\oF}{{\overset{\circ}{\fF}}}

 \def\tb#1{\textcolor{blue}{#1}}
\nc{\fa}{{\mathfrak{a}}}
\nc{\fb}{{\mathfrak{b}}}
\nc{\fg}{{\mathfrak{g}}}
\nc{\fgl}{{\mathfrak{gl}}}
\nc{\fh}{{\mathfrak{h}}}
\nc{\fj}{{\mathfrak{j}}}
\nc{\fm}{{\mathfrak{m}}}
\nc{\fn}{{\mathfrak{n}}}
\nc{\fu}{{\mathfrak{u}}}
\nc{\fp}{{\mathfrak{p}}}
\nc{\fr}{{\mathfrak{r}}}
\nc{\fs}{{\mathfrak{s}}}
\nc{\fsl}{{\mathfrak{sl}}}
\nc{\hsl}{{\widehat{\mathfrak{sl}}}}
\nc{\hgl}{{\widehat{\mathfrak{gl}}}}
\nc{\hg}{{\widehat{\mathfrak{g}}}}
\nc{\chg}{{\widehat{\mathfrak{g}}}{}^\vee}
\nc{\hn}{{\widehat{\mathfrak{n}}}}
\nc{\chn}{{\widehat{\mathfrak{n}}}{}^\vee}

\nc{\fA}{{\mathfrak{A}}}
\nc{\fB}{{\mathfrak{B}}}
\nc{\fD}{{\mathfrak{D}}}
\nc{\fE}{{\mathfrak{E}}}
\nc{\fF}{{\mathfrak{F}}}
\nc{\fG}{{\mathfrak{G}}}
\nc{\fK}{{\mathfrak{K}}}
\nc{\fL}{{\mathfrak{L}}}
\nc{\fM}{{\mathfrak{M}}}
\nc{\fN}{{\mathfrak{N}}}
\nc{\fP}{{\mathfrak{P}}}
\nc{\fU}{{\mathfrak{U}}}
\nc{\fV}{{\mathfrak{V}}}
\nc{\fZ}{{\mathfrak{Z}}}

\nc{\bb}{{\mathbf{b}}}
\nc{\bc}{{\mathbf{c}}}
\nc{\bd}{{\mathbf{d}}}
\nc{\be}{{\mathbf{e}}}
\nc{\bj}{{\mathbf{j}}}
\nc{\bn}{{\mathbf{n}}}
\nc{\bp}{{\mathbf{p}}}
\nc{\bq}{{\mathbf{q}}}
\nc{\bF}{{\mathbf{F}}}
\nc{\bu}{{\mathbf{u}}}
\nc{\bv}{{\mathbf{v}}}
\nc{\bx}{{\mathbf{x}}}
\nc{\bs}{{\mathbf{s}}}
\nc{\by}{{\mathbf{y}}}
\nc{\bw}{{\mathbf{w}}}
\nc{\bA}{{\mathbf{A}}}
\nc{\bK}{{\mathbf{K}}}
\nc{\bI}{{\mathbf{I}}}
\nc{\bB}{{\mathbf{B}}}
\nc{\bG}{{\mathbf{G}}}
\nc{\bC}{{\mathbf{C}}}
\nc{\bD}{{\mathbf{D}}}
\nc{\bP}{{\mathbf{P}}}
\nc{\bH}{{\mathbf{H}}}
\nc{\bM}{{\mathbf{M}}}
\nc{\bN}{{\mathbf{N}}}
\nc{\bV}{{\mathbf{V}}}
\nc{\bU}{{\mathbf{U}}}
\nc{\bL}{{\mathbf{L}}}
\nc{\bT}{{\mathbf{T}}}
\nc{\bW}{{\mathbf{W}}}
\nc{\bX}{{\mathbf{X}}}
\nc{\bY}{{\mathbf{Y}}}
\nc{\bZ}{{\mathbf{Z}}}
\nc{\bS}{{\mathbf{S}}}

\nc{\sA}{{\mathsf{A}}}
\nc{\sB}{{\mathsf{B}}}
\nc{\sC}{{\mathsf{C}}}
\nc{\sD}{{\mathsf{D}}}
\nc{\sF}{{\mathsf{F}}}
\nc{\sG}{{\mathsf{G}}}
\nc{\sK}{{\mathsf{K}}}
\nc{\sM}{{\mathsf{M}}}
\nc{\sO}{{\mathsf{O}}}
\nc{\sQ}{{\mathsf{Q}}}
\nc{\sP}{{\mathsf{P}}}
\nc{\sZ}{{\mathsf{Z}}}
\nc{\sfp}{{\mathsf{p}}}
\nc{\sr}{{\mathsf{r}}}
\nc{\sg}{{\mathsf{g}}}
\nc{\sff}{{\mathsf{f}}}
\nc{\sfb}{{\mathsf{b}}}
\nc{\sfc}{{\mathsf{c}}}
\nc{\sd}{{\ltimes}}

\nc{\tA}{{\widetilde{\mathbf{A}}}}
\nc{\tB}{{\widetilde{\mathcal{B}}}}
\nc{\tG}{{\widetilde{G}}}
\nc{\TM}{{\widetilde{\mathbb{M}}}{}}
\nc{\tO}{{\widetilde{\mathsf{O}}}{}}
\nc{\tU}{\widetilde{U}}
\nc{\TZ}{{\tilde{Z}}}
\nc{\tx}{{\tilde{x}}}
\nc{\tq}{{\tilde{q}}}

\nc{\tfP}{{\widetilde{\mathfrak{P}}}{}}
\nc{\tz}{{\tilde{\zeta}}}
\nc{\tmu}{{\tilde{\mu}}}

 \def\e{\epsilon}

  \nc{\Ob}{{\mathop{\operatorname{\rm Ob}}}}
  \nc{\Sym}{{\mathop{\operatorname{\rm Sym}}}}
   \nc{\Aut}{{\mathop{\operatorname{\rm Aut}}}}
 \nc{\Spec}{{\mathop{\operatorname{\rm Spec}}}}
  \nc{\spec}{{\mathop{\operatorname{\rm Spec}}}}
\nc{\Ker}{{\mathop{\operatorname{\rm Ker}}}}
 \nc{\dom}{{\mathop{\operatorname{\rm dom}}}}
\nc{\End}{{\mathop{\operatorname{\rm End}}}}
 \nc{\Hom}{\on{\Hom}}
 \nc{\GL}{{\mathop{\operatorname{\rm GL}}}}
 \nc{\Id}{{\mathop{\operatorname{\rm Id}}}}
 \nc{\rk}{{\mathop{\operatorname{\rm rk}}}}
 \nc{\length}{{\mathop{\operatorname{\rm length}}}}
\nc{\supp}{{\mathop{\operatorname{\rm supp}}}}
\nc{\val}{{\rm val}}
\nc{\res}{{\mathop{\operatorname{\rm res}}}}

\def\Ind#1#2#3{{#1} {\downarrow}_{#3} {#2} }

\def\tensor{{\otimes}}
\def\meet{\cap}
\def\union{\cup}

\def\si{\sigma}
\def\g{\gamma}
\def\G{\Gamma}

\def\<{\begin}
 \def\>{\end}

\def\m{\smallsetminus}

\nc{\seq}[1]{\stackrel{#1}{\sim}}

\def\inv{^{-1}}
\def\claim#1{{\noindent \bf Claim #1.\ }}

\def\beq#1{{\begin{equation} \label{#1}}  }

\def\Uu{\mathbb U}

\def\prf{\begin{proof}}
\def\pv{\end{proof} }
 \def\eprf{\end{proof} }

\def\acl{\mathop{\rm acl}\nolimits}
 \def\dcl{\mathop{\rm dcl}\nolimits}

\def\liminv{\underset{\longleftarrow}{lim}\,}

 \def\lbl#1{  \label{#1}  }
\def\a{\alpha}

\def\ba{\bar{a}}
\def\k{{\rm k}}
 \renc{\b}{{\beta}}

\def\std#1{{\widehat{#1}}}

\def\tp{\mathrm{tp}}

\def\Ind#1#2{#1\setbox0=\hbox{$#1x$}\kern\wd0\hbox to 0pt{\hss$#1\mid$\hss}
\lower.9\ht0\hbox to 0pt{\hss$#1\smile$\hss}\kern\wd0}

\def\Lam{\Lambda}
\def\L{\Lambda}

 \author{Ehud Hrushovski}
 \title{Imaginaries and definable types in algebraically closed valued fields}
\address{\newline Institute of Mathematics, the Hebrew
 University of Jerusalem, Givat Ram, Jerusalem, 91904, Israel.} 
 
 \email{ehud@math.huji.ac.il}

\begin{document}
\maketitle
This manuscript is largely an exposition of  material from \cite{hhmcrelle}, \cite{hhm} and \cite{HL},
regarding definable types in the model completion of the theory of valued fields, and the   classification of imaginary sorts.  
The proof
of the  latter  is based here on definable types rather than invariant types,
and on the notion of {\em generic reparametrization}; it allows a more conceptual view than we had when  \cite{hhmcrelle} was written.  
I also try to bring out the relation to the geometry of \cite{HL} - stably dominated definable types as the model theoretic incarnation of a Berkovich point.

The text is based on  notes from a class entitled {\em
Model Theory of Berkovich Spaces}, given at the Hebrew University in the fall term of 2009, and retains the flavor of class notes.   Thanks to  Adina Cohen,  Itai Kaplan, and Daniel Lowengrub for comments.  Most recently, Will Johnson went through the notes with great care; he is due
thanks for    numerous textual improvements as well as  some highly perceptive mathematical comments and corrections.   He further
discovered a considerable simplification of the proof of elimination of imaginaries, based on definable types and their coding in $\Oo$-submodules
of finite dimensional $K$-spaces, but shortcutting the decomposition theorem of definable types, \thmref{8.5}; this proof, I hope, will  appear separately.

 The material was discussed in my talk in the Valuation Theory meeting in El Escorial in 2011.  The slides for this talk can be found in \cite{escorial}.

\ssec{Notation}
We will use a universal domain for a given theory,  usually the theory ACVF defined below.
This is a highly saturated and highly homogeneous model, denoted      $\Uu$.  Small subsets of $\Uu$ are denoted by $A$, $B$, \ldots.
   Definable subsets of $\Uu$ are denoted by $X$, $Y$, \ldots, and sometimes $D$.  If $M$ is a model containing the parameters used to define $X$, $X(M)$ denotes the interpretation of $X$ in $M$.  If $A$ is a substructure of a model and $x_1, \ldots, x_n$ are tuples from the model, then $A(x_1,\ldots,x_n)$ denotes the definable closure of $A, x_1, \ldots, x_n$.

When working with valued fields, the valued field itself   is denoted $K$, the residue field is denoted $\k$, the valuation ring is denoted $\Oo$, the maximal ideal is denoted $\Mm$, and the value group is denoted $\G$.  The residue map is $\res : \Oo \to \k$, and the value map is $\val : K \to \G \cup \{\infty\}$.  The value group is written additively, so that $\Oo = \{x \in K : v(x) \ge 0\}$.  ACVF is the theory of non-trivially valued algebraically closed valued fields.

Let $B_n$ denote the group of invertible upper triangular matrices.    The group of elements of $B_n$ with entries in a given ring $R$
is   denoted $B_n(R)$.   We will also write $B_n$ for $B_n(K)$.
$U_n$ is the group of matrices in $B_n$
with $1$'s on the diagonal.  $D_n$ is the group of diagonal matrices, so that $B_n = D_n U_n$.

\section{ Definable types}

\ssec{Definable types}  Let $L_{Y}$ be the set of formulas of $T$ in variables from $Y$, up to $T$-equivalence.
A {\em definable type}  $p(x)$ is a   family of Boolean retraction $L_{x,Y}$ to $L_{Y}$ (for any finite set of variables
$Y$),  compatible with inclusions $Y \subset Y'$.    
 It is denoted:
$\phi \mapsto (d_px)\phi$.  Thus $(d_px)\phi$ is a formula  with (at most) the same $y$-variables but without the free variable $x$;
it is analogous to quantifiers, but simpler; one says: for generic $x \models p$, $\phi$ holds.

Given a definable type $p$ and a substructure $A$ of $M \models T$, we let
$$p | A = \{\phi(x,a):  a \in A^l, M \models (d_px) \phi (a) \}$$

So we can think of a definable type as a compatible family of types, given systematically over all base sets.

\ssec{Examples, notation} While the development is at first abstract, we will give examples from ACVF, the theory of algebraically
closed valued fields.  $K$ denotes the field, $\Oo$ the valuation ring, $\G$ the value group, $\val$ the valuation map, $\res$ the residue
map into the residue field $k$.

\ssec{Pushforward of definable types}

Let $f: X \to Y$ be an $A$- definable function, and $p$ an $A$- definable type on $X$.  Define $q=f_*p$ by:
$$(d_q y) \theta(y,u) = (d_p x) \theta(f(x),u)$$

Excercise.    For any $B$ containing $A$ we have:  $(f_*p)|B =tp(f(c)/B)$ where $c \models p|B$.

\ssec{Product of definable types} If $p$ and $q$ are two $A$-definable types, then the product $p(x) \otimes q(y)$ is defined by
$$ (d_{p \otimes q} (x,y)) \theta(x,y,u) = (d_q y)(d_p x)\theta(x,y,u). $$
If $B$ contains $A$, then $(c_1,c_2) \models p \otimes q | B$ if and only if $c_2 \models q | B$ and $c_1 \models p | B(c_2)$.

\ssec{Orthogonality}   A definable type $q(x)$ is {\em constant} if $(d_q x)(x=y)$ has a solution.

Excercise.  In this case, $(d_q x)(x=y)$ has a unique solution $a$; and $a$ is the unique realization of
$q|B$, for any $B$ over which $q$ is defined.

\<{defn}
$p$ is {\em orthogonal to $\G$} if for any $\Uu$-definable function $f$ into $\G$,
 $f_*p$ is constant.  \>{defn}

Equivalently, by considering coordinate projections, any $\Uu$-definable function $f$ into $\G^n$ is constant.
   We will use this definition for the value group, which eliminates imaginaries; otherwise
 we would instead consider  definable functions $f$ into $\G^{eq}$.

\ssec{Stable embeddedness} \label{stable-embeddedness}

A sort $D$ is   {\em stably embedded} if any $\Uu$-definable subset of $D^m$ is $D(\Uu)$-definable.

In ACVF, both $\G$ and $\k$ are {\em stably embedded}; this is an immediate consequence of quantifier-elimination
in the standard three-sorted language (See Theorem 2.1.1 (iii) in \cite{hhmcrelle}, or the first paragraph of the Appendix.)
It suffices to consider atomic formulas, with some variables from $\G$ and some from other sorts.
Any atomic formula $\phi(x_1,\ldots,x_n,y_1,\ldots,y_n)$ with $x_i$ in $\G$, $y_j \in VF$,
has the form:  $\theta(x_1,\ldots,x_n,\val(h_{\nu}(y)))$.  So $\phi(x,b)$ defines the same set as
$\theta(x,d)$ where $d = \val h(b)$.  Similarly for $\k$ and $\res h$, with $h$ a rational function.

    Orthogonality of $p$ to $\G$ can also be stated as follows:  Let $B'=B(\g)$ be generated over  $B$ by some realizations of $\G$.
Then $p|B$ implies $p|B'$.   

\ssec{Domination}

\<{lem}  \lbl{dom1}  Let $f: X \to Y$ be an $A$-definable function.  Let $q$ be an $A$-definable type on $Y$, and let $p_A$ be a type over $A$ on $X$.  
Assume:     for any $B \geq A$ there exists a unique type $p_B$  such that
$p_B$ contains $p_A$, and $f_* p_B = q|B$.   Then there exists a unique $A$-definable type $p$ such that for all $B$, $p|B=p_B$.  \>{lem}

\prf  More generally, let us say a type  $p_\Uu$ over $\Uu$ is {\em definably generated } over $A$  if it is generated by a    partial  type of the form $\union_{(\phi,\theta) \in S} P(\phi,\theta)$, where $S$ is a (small) set of pairs of formulas
$(\phi(x,y), \theta(y))$ over $A$, and $P(\phi,\theta) = \{\phi(x,b): \theta(b) \}$.

It sufices to show that if $p_\Uu$ is definably generated over $A$, then $p_\Uu$ is definable over $A$, i.e.
$\{b: \phi(x,b) \in p_\Uu \}$ is $A$-definable for each $A$-formula $\phi(x,y)$.

Let $\phi(x,y)$ be any formula.   From the
fact that $p_\Uu$ is definably generated it follows easily that
$\{b: \phi(x,b) \in p_\Uu \}$ is an $\bigvee$-definable set over $A$, i.e.  a union of $A$-definable sets.   Indeed,  $\phi(x,b) \in p_\Uu$ if and only if for some $(\phi_1,\theta_1),\ldots,(\phi_m,\theta_m) \in S$, $(\exists c_1,\cdots,c_m)(\theta_i(c_i) \wedge (\forall x) (\bigwedge_i \phi_i(x,c) \implies \phi(x,b) )$.
Applying this to $\neg \phi$, we see that the complement of $\{b: \phi(x,b) \in p_\Uu \}$ is also $\bigvee$-definable.
Hence $\{b: \phi(x,b) \in p_\Uu \}$ is $A$- definable.

 \eprf

\<{defn}  In the situation of the lemma,  $p$ is said to be {\em dominated by} $q$ via $f$ \>{defn}

In the situation of the lemma, $p$ is said to be {\em dominated by} $q$ via $f$.  More precisely: 
\<{defn}    $p$ is dominated by $q$ via $f$ if there is some $A$ over which $p$, $q$, and $f$ are defined, such that for every $B \ge A$,
$(q|B)(f(x)) \cup (p|A)(x) \vdash (p|B)(x)$. \>{defn}

In general, when $p,q,f$ are $A$-definable, one can visualize that $p$ is dominated by $q$ over
some bigger set $B$, but not over $A$.  When $A$ is a model, this does not happen, nor will it occur
in our setting of stable domination (see Remark \ref{3.10}).   (Thanks to Will Johnson for this remark.)

\<{example}  (ACVF)  Let $X=\Oo, Y = \k, f = \res$.   Let $q$ be the generic type of $\k$, i.e. $q|B$
is generated by:  $y \in \k, y \notin V$ for any finite $B$-definable $V$.   Then $x \in \Oo,  f(x) \models q|B$ generates
a complete type $p|B$ over $B$.   This is called the {\em generic type of $\Oo$}.
\>{example}

\<{ex}\lbl{ball-generic}  Show that  $p|B$ is complete.  For any polynomial $\sum b_i x^i$ over $B$, show that
$\val( \sum b_ix^i) = \min_i  \val (b_i)$ \tb{for $x$ realizing $p | B$.  In particular, $p$ is orthogonal to $\G$}.\>{ex}

\<{example}  Let $\Mm=\{x: \val(x)>0 \}$ be the maximal ideal.  Let $f(x)=\val(x)$.  Let $q(x)$ be the type just
above $0$ in $\G$.  Then $q$ dominates via $f$ a definable type $p_\Mm$, the generic type of $\Mm$.  \>{example}

\<{example}  (ACVF$_{0,0}$).    Let $a_0,a_1,\ldots \in \Qq$.  Let $\val(t)>0$.
Let $p_0(x,y)$ consist of all formulas (over $\Qq(t)$)

$$ \val( y -  \sum_{k=0}^n  a_k (xt)^k  \geq n \val(t)) $$

Then $p_0(x,y) + (p_\Oo |\Uu) (x)$ generates a complete type $p|\Uu$, provided $\sum a_n x^n$ is transcendental.

Let $X = \Oo \times \Oo$, $Y=\k$, $f(x,y) = \res(x)$.  Then $p$ is dominated by the generic type of $\k$, via $f$.

\>{example}

To prove the domination, say  $\val(t)=1$.     First let $M$ be a   valued field   extension of $\Qq(t)^{alg}$
such that $\Zz$ is cofinal in $\val(M)$.
  We   prove domination over $M$.

 Generalizing the construction, allow $a_n \in \Oo_M$, $a= \sum a_k X^k$, and define $p_0^a$
to consist of all formulas:
$$ \val( y -  \sum_{k=0}^n  a_k (xt)^k ) \geq n $$
For $a$ fixed, write $p_0=p_0^a$.  

Let $c \models p_\Oo | M$.  First suppose $p_0(c,0)$ holds.  Then $\min_{i \leq n} \val(a_i) + i = \val(    \sum_{k=0}^n  a_k (ct)^k)  \geq n  $.
So $\val(a_i) \geq n-i$.  Letting $n \to \infty$ (and using $a_i \in M$) we see that $a_i=0$; so $a=0$.

Next suppose just that $p_0(c,d)$ holds for some $d \in M(c)^{alg}$.  So $F(c,d)=0$ for some polynomial $F \in \Oo_M[x,y]$.
Let $a'=F(x,a(x))$ be the power series obtained by substituting $a(x)$ for $y$.  Let $p_0' = p_0^{a'}$.  Then $p_0'(c,0)$ holds.
Hence by the previous paragraph, $a'=0$, so $a$ is algebraic.

Otherwise, $p_0(c,y)$ defines an infinite intersection $b$ of balls over $M(c)$, with no algebraic point.  Hence $b$ contains no nonempty $M(c)$-definable subset ($M(c)^{alg} \models ACVF$, so any nonempty $M(c)$-definable set does have an algebraic point.)    So $p_0 + tp(c/M)$ generates a complete type over $M(c)$, as promised.

We can take $M$ to be maximally complete; this suffices to    show that  $p|M$ is stably dominated.

Now if $N$ is a valued field extension of $M$ with $\res(N)=\res(M)$, then $p|M \vdash p|N$, hence $p_0(x,y) + p_\Oo |M$ already
generates $p|N$.

But any valued field extension of $\Qq(t)^{alg}$ can be obtained in this way (taking such an $M,N$ and then a  subextension.)  This proves
the domination statement in the example.

\ssec{Density of definable types}

We consider the following extension property for a definable set $D$
over a base set $A$, possibly including imaginaries.

 (E(A,D)): Either $D = \emptyset$, or there exists a definable type $p$ on $D$ (over $\Uu$) such that
$p$ has a finite orbit under $Aut(\Uu/A)$.

Say $T$ has property $E$ if $E(A,D)$ holds for all $A,D$.   In \lemref{8.1} below, we will see that ACVF has property (E).

 
We say that a substructure $B$ of $\Uu$ is a {\em canonical base} for an object $p$ constructed from $\Uu$  if for any $\si \in Aut(\Uu)$, $\si(p)=p$ iff $\si(b)=b $ for all $b \in B$.

\<{lem}\lbl{dd3} Let $T$ be a theory with property (E), and assume any definable type (in the basic sorts) has a canonical base in certain imaginary sorts $S_1,S_2, \ldots$.
Then $T$ admits elimination of imaginaries to the level of finite subsets of products of the $S_i$. \>{lem}


\ssec{Definable types on $\G^n$}

Let $\G$ be a divisible ordered Abelian group.  Recall that the theory of divisible ordered Abelian groups has quantifier-elimination
(a result whose roots go back to Fourier.)

We will  consider projections $\phi^a: \G^n \to \G$,
$\phi^a(x) = a \cdot x$, where $a \in \Qq^n \m (0)$.

We say two definable types $p,q$ are {\em orthogonal} if 
 there is a set $A$ over which $p$ and $q$ are defined, such that for
any $B \geq A$,
$p(x) |B \union q(y) | B$ generates a complete type in  the  variables $x,y$.

 A definable type $p$ in $\G^n$ has a {\em limit} if there is some $c \in \G^n$ such that for every $\Uu$-definable open neighborhood $U$ of $c$,
the formula $x \in U$ is in $p | \Uu$.

\<{lem} \lbl{cc5} Let $p$ be a definable type of $\G^n$, over $A$.  Then up to a change of coordinates by
a rational   $n \times n$ matrix,
  $p$ decomposes as the join of two orthogonal definable types $p_f, p_i$, such that $p_f$ has a limit in $\G^m$,
and $\phi^a_* p_i$ has limit point $\pm \infty$ for any $a \in \Qq^n  \m (0)$. \>{lem}

\prf   Let
$\a_1,\ldots,\a_k$ be a maximal set of linearly independent vectors in $\Qq^n$ such that
the image of $p$ under $(x_1,\ldots,x_n) \mapsto   \a_i \cdot x$ has a limit point in $\G$ 
\footnote{Equivalently, the image of $p$ under $x \mapsto (\alpha_1x, \ldots, \alpha_k x)$ 
has a limit point in $\G^k$.}
Let $\b_1,\ldots,\b_l$ be a maximal set of vectors in $\Qq^n$ such that for {any/every model $M$ and for} $ x \models p|M$,
$\a_1 x ,\ldots,\a_k x ,\b_1 x,\cdots,\b_l x$ are linearly independent over  $M(\a_1 x, \ldots, \a_k x)$
 If $a \models p|M$, let $a'=(\a_1 a, \ldots, \a_k a)$, $a'' = (\b_1 a, \ldots,\b_k a)$.    For $\a \in \Qq(\a_1,\ldots,\a_k)$, 
the element
 $\a a$ is bounded between elements of $M$.  On the other hand each $\b a$ ($\b \in \Qq(\b_1,\ldots,\b_k)\m \{0\})$ satisfies
 $\b a > M$ or $\b a < M$.
For if $m \leq \b a'' \leq m'$ for some $m \in M$, since $\tp(\b a''/M)$ is definable it must have a finite limit,
contradicting the maximality of $k$.  It follows that   $\tp(\a a/  M) \union \tp(\b a /M)$ extends to a complete 2-type, namely $\tp((\a a ,\b a) / M)$; in particular $\tp(\a a + \b a / M)$
is determined; from this, by quantifier elimination, $\tp(a'/M)  \union \tp(a''/M)$ extends to a unique type in $k+l$ variables.
So $\tp(a'/M)$, $\tp(a''/M)$ are orthogonal.  After some sign changes in $a''$, so that each coordinate is $>M$,
the lemma follows.
\eprf

\<{lem} \lbl{cc6}  \<{enumerate} \item  Let $p,p'$ be   definable types on $\G^n$. If $\phi^a_*p=\phi^a_*p'$
for each $a$, then $p=p'$.
\item
Let $p$ be a definable type on $\G^n$.  If $\phi^a_*p$ is 0-definable for each $a$, then $p$
is $0$-definable.  \>{enumerate} \>{lem}

\prf  (1) Any formula $\phi(x,y)$ is a Boolean combination of formulas $a \cdot x + b \cdot y > c$ (or $=c$).
The definition of such a formula is determined by $\phi^a_* p$.

 (2) Let $\si$ be an automorphism, $p'=\si(p)$; we have to show that $p'=p$.  This follows from (1).
\eprf

\<{lem} \lbl{cc7} Let $p$ be a definable type of $\G^n$.  For $c \in \G^n$, let  $\a^c(x)=x {+} c$ .  Then for some $c \in \G^n$, $\a^c _* p$ is $0$-definable.
\>{lem}

\prf A linear change of coordinates (with $\Qq$-coefficients) does not effect this statement.  So
we may assume the conclusion of \lemref{cc5} holds.   Translating the $p_f$ part by $-\lim p_f$, we may assume
$p_f$ has limit $0 \in \G^m$.  It follows that for any $a \in \Qq^n \m (0)$, $\phi^a_* p$ has limit $0$ or $\pm \infty$.
There are only five definable 1-types with this property, all 0-definable.  Hence by \lemref{cc6}(2), $p$ is 0-definable.  \eprf

\section{Algebraic lemmas on valued fields}

The material in this section is classical, going back in part to Ostrowsky and Kaplansky; see the book by F.V.-Kuhlmann {\tt{http://math.usask.ca/~fvk/Fvkbook.htm}.}

\<{defn}  An extension $L \leq L'$ of valued fields is {\em immediate} if $L,L'$ have the same value group and residue field.

$K$ is {\em maximally complete} if it has no proper immediate extensions.  \>{defn}

\<{ex} \lbl{a1} Let $K$ be an algebraically closed valued field,  $L$ a valued field extension, $t \in L$.  Assume
$L=K(t)$ as a field.    Since any element
of $K[t]$ is a product of linear factors, the valuation on $L$  is determined by
$v(t-a)$ for $a \in K$.   Then one of the following holds:  \<{itemize}
\item  $v(t-a) = \g \notin \G(K)$ for some $a \in K$.    Show that $\G(L) = \G(K)(\g)$, $\k(L)=\k(K)$.
\item  $v(t-a)\in \G(K)$ for all $a \in K$, and $v(t-a)$ takes a maximal value $v(b)$ at some $a \in K$.  Show that
$\k(L) = \k(K)(e)$ where $e = \res ((t-a)/b)$.
\item $v(t-a)\in \G(K)$ for all $a \in K$, and a maximum is not attained.  Show that   $K(t)$ is an immediate extension.
\>{itemize}  \>{ex}

\<{lem} \lbl{a2}  Let $L/K$ be an extension of valued fields.  Then $tr.deg._{\res K} \res L + \dim_\Qq (\val(L)/\val(K)) \leq tr. deg._K(L)$.  \>{lem}

\prf  This reduces to the case that $L/K$ is generated by one element.  In this case $L/K$ is algebraic or $L=K(t)$ is a rational function field.  In the algebraic case, $\res L $ is a finite extension of $\res K$ (of some degree $e$)
and $\val(L) /\val(F)$ is finite (of some order $f$; in fact we have $ef \leq [L:K]$.)    In case $L=K(t)$, we may assume $K$ is algebraically closed, since passing to this case will not lower the left hand side; and Ex. \ref{a1} applies.
\eprf

\<{lem}   \lbl{a3}  Let $K$ denote a valued field, with algebraically closed residue field $\k$ and divisible value group $A$.  Assume $K$ is maximally complete,    \<{itemize}
\item    $K$ is algebraically closed.
 \item $K$ is {\em spherically complete}, i.e. any set of balls, linearly ordered by inclusion, has nonempty intersection.

 \>{itemize}  \>{lem}
 \prf  (1)  This follows from \lemref{a2}:  algebraic extensions are immediate since the value group and residue field have no proper finite extensions.  (2)  Let $b_i$ be a set of balls, indexed by a linear ordering $I$.  If $\meet_i b_i = \emptyset$, then for any $a \in K$ we have $a \notin b_i$ for large $i$, and it follows that $\alpha(a) = v(a-c)$ is constant for $c \in b_i$.    Define a valuation on $K(t)$
 by $v(t-a) = \alpha(a)$.    Then by Ex. \ref{a1}   this
   is an immediate extension, a contradiction.  \eprf

Any valued field $K$ has a maximally complete immediate extension, of cardinality at most $2^{|K|}$.

\ssec{Valued vector spaces}
A {\em valued} vector space over  valued field $K$ is a triple $(V,\G(V),v)$, with $V$ a $K$-space,
$\G(V)$ a  linearly ordered set $\G(V)$ with an action  $+: \G(K) \times \G(V) \to \G(V)$, order-preserving in each variable, and $v$ a map
 $v: V \m (0) \to \G(V)$ with $v(a+b) \geq \min (v(a),v(b))$ and $v(cb) = v(c)+v(b)$ for $a,b \in V, c \in K$.

 If $a_1,\ldots,a_n$ are elements of $V$ with $v(a_1),\ldots,v(a_n)$ in  distinct $\G(K)$-orbits,
 it follows that $a_1,\ldots,a_n$ are linearly independent over $K$.  In particular if $V$
  is finite-dimensional, $\G(V)$ can only consist of
 finitely many $\G(K)$-orbits.

By a {\em ball} in $V$ we mean a set of the form $\{b \in V: v(a-b) \geq \alpha \}$.
$V$ is {\em spherically complete} if any set of balls, linearly ordered by inclusion, has nonempty intersection.

A set $a_1,\ldots,a_n$ of elements of $V$ is   called {\em separated}
if  for all $c_1,\ldots,c_n \in K$,   we have
$$v  \left( \sum c_i a_i \right)   =  \min_i ( v(c_i)+v(a_i))$$
Such a set is in particular linearly independent.

If $V=K^n$ is a valued $K$-space with a  separated basis, a ball for $V$ is just a product of balls of $K$,
so $V$ is spherically complete if $K$ is.

 {If $V$ is a valued $K$-space with a spherically complete subspace $W \le V$, and $a \in V$, then the set $\{ v(w - a) : w \in W\}$
attains a maximum, because for each $\gamma \in \Gamma(V)$, the set $\{w \in W : v(w - a) \ge \gamma\}$ is either empty or a
ball in $W$.}

\<{lem}  \lbl{a4} Let $K$ be a spherically complete valued field, $V$ a finite-dimensional $K$-space.  Then $V$
   has a separated basis.     \>{lem}

\prf

Let $a_1,\ldots,a_m$ be a maximal separated set, $U$ the subspace generated by  $a_1,\ldots,a_m$.
Then $U$ has a separated basis, so it is spherically complete.
 If $U =V$ we are done.  Otherwise, let $a \in V \m U$.     Consider  the possible values $v(u-a)$, $u \in U$.   Since
$U$ it is spherically complete, so there must
be a maximal value among these.    Replacing $a$ by $a-u$ with $v(a-u)$ maximal, we may assume $v(a) \geq v(a-u)$ for all $u \in U$.  In this case, $a_1,\ldots,a_m,a$ is separated.  For given $c_1,\ldots,c_m$, we have
$v(\sum c_i a_i) = \min_i v(c_ia_i) = \gamma $ say.
It suffices to see that $v(\sum c_i a_i + a) \leq \min (\gamma,v(a))$;   this follows from the strong triangle inequality
 when $\gamma \neq v(a)$, and  from $v(a) \geq  v(\sum c_i a_i + a) $ when $\g=v(a)$.   \eprf

\ssec{Induced $\k$-spaces}

Let $V$ be a valued $K$-space, and $\a \in \G(V)$.  Then $\Lambda_\a = \{a \in V: v(a) \geq \a \}$ is an $\Oo$-submodule,
and $\Lambda_\a^o = a \in V: v(a) > \a \}$ is an $\Oo$-submodule containing $\Mm \Lambda_\a$.
Let $V_\a = \Lambda_\a /\Lambda_\a^o$; this is a $\k=\Oo/\Mm$-space, finite-dimensional if $V$ is.

Let $h: U \to V$ be a homomorphism of valued $K$-spaces; meaning there is also a map $h: \G(U) \to \G(V)$
of $\G(K)$-sets, with $h(\a) < h(\b)$ when $\a < \b$, and $v(h(a)) = h(v(a))$.  Then
$h$ induces a homomorphism $U_\a \to V_{h(\a)}$ for each $\a$.

\ssec{Tensor products}

    Let $U,V$ be valued $K$-spaces.   Consider
$K$-spaces $(W,\G(W))$ and maps
$$h: U \tensor V \to W, \ \  +: \G(U) \times \G(V) \to \G(W)$$
such that $v(h(a \tensor b)) =v(a)+v(b)$ and $\G(U) \times \G(V) \to \G(W)$ is order-preserving in each variable.

 Then for each $\a \in \G(U), \beta \in \G(V)$ we have an induced homomorphism
 $$U_\a \tensor V_\b \to W_{\a + \b}$$

 \<{lem}  \lbl{a5} Let $K$ be spherically complete, and let $U,V$ be valued $K$-spaces.  Let $E$ be a 
 divisible ordered Abelian group with $\G(K)$-action, and assume $\G(U) , \G(V) \leq E$ and 
 $E \cong \G(U) \times_{\G(K)} \G(V)$, i.e. if $\a+\b=\a'+\b'$ then for some $\g \in \G(K)$, $\g+\a=\a'$
and $\g+\b'=\b$.  Then 
 there exists a unique $(W,h:W \to E)$ (up to a unique isomorphism) such that:
   \<{enumerate}
\item  For any $(\a,\b) \in \G(U) \times \G(V)$, the induced homomorphism $U_\a \tensor V_\b \to W_{\a + \b}$
is injective.
\>{enumerate}
  \>{lem}

\prf  To prove uniqueness
we have to show that $h$ is injective, and determine $v(h(x))$ for all $x \in U \tensor V$.
Write $x = \sum_{i=1}^n a_i \tensor b_i$ where $(a_1,\ldots,a_n)$ are separated.  Then it suffices to show:

\claim{}  $h(x) \neq 0$, and  $v(h(x)) = \min_i v(a_i) + v(b_i)$.

By grouping the terms according to the value of $v(a_i) + v(b_i)$, it suffices to prove the claim when $v(a_i)+v(b_i)$ is constant.
In this case by assumption there exists $\g_i=v(c_i), c_i \in K$ with $v(a_i)=v(a_1)+\g, v(b_i)=v(b_1)-\g$.
Replacing $a_i$ by $a_i/c_i$ and $b_i$ by $b_ic_i$, we may assume $v(a_i)=v(a_1)=\a, v(b_i)=v(b_i)=\b$.
So $a_i \in \Lambda^U_\a, b_i \in \Lambda^V_\b$.  Since $a_1,\ldots,a_n$ are separated, the images $\ba_i$ of the $a_i$ in    $U_\alpha$  are linearly independent.    The images $\bar{b}_i$ of the $b_i$ in $V_\b$ are nonzero.  Hence $\sum \ba_i \tensor \bar{b}_i \neq 0 \in U_\a \tensor V_\b$.  Since $h$ induces an injective map into $W_{\a+\b}$  it follows that
 $v(h(\sum_i a_i \tensor b_i)) = \a+\b$.

 With uniqueness proved, functoriality is clear and so it suffices to prove existence in the finite dimensional
 case.    This is easily done by choosing
a separated basis and following the recipe implicit above.
  \eprf

\<{prop} \lbl{a6} Let $K$ be a spherically complete valued field, $L_1,L_2$ valued field extensions,
within a valued field extension $N$ generated by $L_1 \union L_2$.
Assume $\G(K)=\G(L_1)$, and $\k(L_1)$ is linearly disjoint from $\k(L_2)$ over $\k(K)$.
Then the structure of the valued field $N$ is uniquely determined given $L_1$ and $L_2$.
\>{prop}

\prf    It suffices to show that the natural map $h:L_1 \tensor L_2 \to N$ is injective and that
$v(h(x))$ is determined for $x \in L_1 \tensor L_2$, since passage to the field of fractions is clear
using $v(x/y) = v(x)-v(y)$.  Let $W$ be the image of $h$.  Then we are in the setting of \lemref{a5}, (1) holds,
and (2) is clear since $\G(U)=\G(K)$.  For the same reason, (3) reduces to the case $\a=0$.
Suppose $\val(b_1) = \min_i \val(b_i)$, without loss of generality.
We have $h(\sum_i a_i \tensor b_i) = b_1 h( \sum_i a_i \tensor (b_i/b_1))$, so we may take $\b=0$ too.
In this case (3) amounts to the linear disjointess assumption.  The corollary now follows from the lemma.
\eprf

\propref{a6} will imply that any definable type orthogonal to  $\G$ is dominated by its images in $\k$.
We did not use \lemref{a5} in full generality; using it we could deduce
that any definable type is dominated by its images in $\G$ and in $\k$.    We will in fact require
a stronger statement, of stable domination relative to $\G$.  The algebraic content consists of the lemma below.

Let  $L_1,L_2$ be two  valued field extensions $L_1,L_2$ of
a valued field $K$, contained in a valued field extension $N$ of $K$, and such that
$L_1 \union L_2$ generates $N$.
As in \lemref{a6}, we will say that the interaction between  $L_1,L_2$ is  uniquely determined
(given some conditions)
 if whenever $N'$ is another valued field extensions of $K$, and
$j_i: L_i \to N'$   are valued $K$-algebra homomorphisms (satisfying the same conditions),
then   there exists a (unique)  valued $K$-algebra embedding $j: N \to N'$ with $j |L_i = j_i $.

It is easy to see that condition (2) below does not depend on the choice of $Z$.

\<{prop} \lbl{a8} Let $K$ be a spherically complete valued field, $L_1,L_2$ valued field extensions,
within a valued field extension $N$ generated by $L_1 \union L_2$.    Let $\k_0=\k(K)$, $\k_i = \k(L_i)$, $\k_{12} =\k_1\k_2$.
Then the interaction of $L_1,L_2$ is uniquely determined assuming the following conditions.  \<{enumerate}
\item  $\G(L_1) \subseteq \G(L_2)$.
\item Let $Z$ be a $\Qq$-basis for  $\G(L_1) / \G(K)$;
for $z \in Z$ let $a_z \in L_1$ and $b_z \in L_2$ have $v(a_z)=v(b_z)=z$, and let $c_z=a_z/b_z$.
Assume the elements $\res(c_z)$ form an algebraically independent set over $\k_{12}$.
\item $\k_1,\k_2$ are linearly disjoint over $\k_0$.
\>{enumerate}
\>{prop}

\prf    As in \lemref{a6}, it suffices to show that the natural map $h:L_1 \tensor L_2 \to N$ is injective and that
$v(h(x))$ is determined for $x \in L_1 \tensor L_2$.   Write $x=\sum_{i=1}^n a_i \tensor b_i$, with $(a_i)$ separated.   We claim
that $v(x) = \min_i v(a_i)+v(b_i) \in \G(L_2)$.  As before we may assume $v(a_i)+v(b_i)=\g$ does not depend on $i$.
Moreover since $\g=v(c)$ for some $c \in L_2$, dividing $b_i$ by $c$ we may assume $\g=0$.  The subgroup of
$\G(L_2)$ generated by the $v(a_i)$ is finitely generated; let $d_1,\ldots,d_l$ be a minimal set of generators
of this group modulo $\G(K)$.  Let $a'_z,b'_z,c_z$ be as in condition (2), so that $v(a'_z)=v(b'_z)=d_z$ for
$z=1,\ldots,l$, $c_z=a'_z/b'_z$, and the elements $\res(c_z)$ are algebraically independent over $\k_{12}$.
For each $i$, there exists $m=(m_1,\ldots,m_l) \in \Zz^l$ with $ \sum m_zd_z = v(a_i)$.  Write $m(i)$ for this $m$,
and $(a')^{m}$ for $\Pi_z (a'_z)^{m_z}$, and similarly for $b'$ and $c$.   Let $A_i = a_i / (a')^{m(i)},
B_i = b_i (b')^{m(i)}$.  Then
$$\sum a_ib_i = \sum A_iB_i c^{-m(i)} $$
We have to show that this has valuation zero, i.e. that \[ \sum_i \res(A_i) \res(B_i) \res(c)^{-m(i)}  \neq 0 \]
Since the $\res(c_z)$ are algebraically independent (2), it suffices to show that for a fixed value of $m \in \Zz^l$, we have:
$\sum_{m(i)=m} \res(A_i)  \res(B_i) \neq 0$.  But this follows from (3) as in \lemref{a5}.  \eprf

\section{Stably dominated types}.

\<{defn}  An $A$-definable type $p$ is stably dominated if for some $B \geq A$, $p$ is dominated over $B$ by
a definable map $f$ into $V$ for some finite-dimensional $\k$-space $V$.  \>{defn}

When the base $A$ consists of elements of the valued field and $\Gamma$, it can be shown
that $f$ can be chosen to be $A$-definable.  The space $V$ is isomorphic to $k^m$ over some larger $B$,
but not necessarily over $A$.  For instance, given $\a \in \G$, let $\Oo \a = \Oo c$ where $\val(c)=\alpha$.
Then $\Oo \a $ is a free $\Oo$-module, and $\Oo a / \Mm a$ is a one-dimensional $\k$-space $V_\a$.

\<{ex}  The generic type of the ball $\val(x) \geq \a$ is dominated by the map $r_\a: \Oo a \to V_\a$.    However
every $\a$-definable map on $\Oo_\a$ into $\k$ is constant, if $\a$ is not a root of the valuation of some element
of the prime field.
\>{ex}

This special family of definable types will be the main object we will look at.
For any definable set $V$, we will define $\std{V}$ to be the set of stably dominated types on $V$.  Later,
a topology will be defined on $\std{V}$; $V$ will be dense in $\std{V}$, called the {\em stable completion} of $V$.

\<{thm} \lbl{b1} In ACVF, the following conditions on a definable type are equivalent:
\<{enumerate}
\item  $p$ is stably dominated.
\item For all definable $q$,  $p(x) \tensor q(y) = q(y) \tensor p(x)$
\item  $p$ is symmetric:  $p(x) \tensor p(y) = p(y) \tensor p(x)$.
\item   $p$ is orthogonal to $\G$.
\>{enumerate}
\>{thm}

\prf   (1) implies (2):   by domination it suffices to prove that $p(x) \tensor q(y) = q(y) \tensor p(x)$
for $p$ on $\k^n$.  By stable embeddedness one reduces to the case that $q$ too is on $\k^n$.

 (2) implies (3) is trivial.

(3) implies (4):  Let $f$ be a definable function into $\G$.  Then $q=f_* p $ is symmetric.  But by considering
the $q(u)$-definition of $u<v$ one sees that $q$ must be  constant.

(4) implies (1).  Let $M$ be a maximally complete valued field, with $p$ definable over $M$.  Let $a \models p|M$,
$N=M(a)$.
Then $\G(N)=\G(M)$ by orthogonality.  By \propref{b3}, a unique $M$-definable type
extends $p|M$, and this type is stably dominated; this type must be $p$.   \eprf

\<{ex} \lbl{b2} Let $k$ be an algebraically closed field, $V$ a finite-dimensional vector space over $k$, definable in some theory over a base $A$.  We assume that the definable subsets of $k^m$ are the constructible subsets.
Let $p_A$ be a type of elements of $V$, over $A$.   Then there exists at most one $A$-definable type $p$
such that $p|A =p_A$.  \>{ex}

(Proof:  $p$ is the generic type of a unique Zariski-closed subset $W$ of $V$; $W$ must be $A$-definable;
we must have $W \in p_A$ but no smaller subvariety is in $p_A$; this characterizes $W$ and hence $p$.)

\<{prop} \lbl{b3}  Let $M$ be a maximally complete algebraically closed valued field, $N=M(a)$ a valued field extension.
Let $\g=(\g_1,\ldots,\g_n)$ be a basis for $\G(N)/\G(M)$.
Then there exists a unique $M(\g)$-definable  type  extending $tp(a/M(\g))$.  This type is stably dominated.
  \>{prop}

\prf  We have $\g_i = \val(c_i)$ for some $c_i \in N$.
Let $e_i = r_{\g_i}(c_i) \in V_{\a_i}$ (see notation above.)
Let $\a_1,\ldots,\a_m$ be a transcendence basis for
$\k(N)$ over $\k(M)$.
We have $\e_i=g_i(a), \a_j = h_j(a)$ for some $M(\g)$-definable functions $g_i,h_j$.  Let
$V=\Pi_i V_{\a_i} \times \k^m$,  $g=(g_1,\ldots,h_m)$.      
Let $q=q(v_1,\ldots,v_n,t_1,\ldots,t_m)$ be the generic type of the $k$-space $V$; equivalently, letting $
V_{i+n} = \k$,
$q = q_1 \tensor \ldots \tensor q_{n+m}$ where $q_i$ is the unique non-constant definable type on the
1-dimensional $k$-vector space $V_i$.   Note that for any  structure $B \geq M(\g)$, for $i \leq n$,
if $e_i ' \in V_i(B), e_i' \neq 0$ and $e_i \models q_i | B$ then $e_i / e_i ' $ is a well-defined element, realizing the generic
type of $\k$ over $B$; hence if   $(e_1,\ldots,e_n,\a_1,\ldots,\a_m) \models q|B$, then
$(e_1/e_1',\ldots,e_n/e_n', \a_1,\ldots,\a_m) $ are algebraically independent over $B$.
Note also in this situation that  if  if $\g_i = \val(b_i)$ with $b_i \in B$,   then $e_i/e_i' = \res(c_i/b_i)$.
  By \propref{a8}, there exists a unique type $p_B$ extending $tp(a/M(\g))$ and with $g_*p_B = q|B$.    By \lemref{dom1}
  there exists a unique $M(\g)$-definable type $p$ with $p|B=p_B$ for all $B$.

  By definition, $p$ is dominated by $g$ and hence stably dominated.  If $p$ is another $M(\g)$-definable type extending
  $tp(a/M(\g))$, let $q'=g_* p'$.  Then $q'$ is an $M(\g)$-definable type extending $tp(e_1,\ldots,\a_m)/M(\g)$.  By
  Excercise \ref{b2} we have $q'=q$, and hence by the domination, $p'=p$.  This proves the uniqueness of $p$.
\eprf

\<{discussion} \rm
Let $V$ be an $M$-definable set, with $a \in V$.
We will see below that $\std{V}$ can be viewed as a  {\em   pro-definable set};  i.e. an inverse limit of definable sets.  In more detail:  
we will describe
certain    definable sets $\std{V}_d$ for $d \in \Nn$,  and definable maps $\std{V}_{d+1} \to \std{V}_d$.  (These maps can be 
taken to be surjective, but we will not use this fact here.  The $\std{V}_d$ will be subsets of $K^m \times S_n$ for appropriate $m,n$,
where $S_n$ is the sort of lattices in $K^n$, described below.)   Let $ \liminv_{d \in \Nn} \std{V}_d$ be the set of sequences
$c=(c_d: d\in \Nn)$ such that $c_{d+1} \mapsto c_d$.  Say  $c \in \dcl(A)$ iff each $c_d \in \dcl(A)$.  
A definable map $f: X \to  \liminv_{d \in \Nn} \std{V}_d$ means:   a compatible system of definable maps $f_d: \G^n \to \std{V}_d$.

For each  $c \in   \liminv_{d \in \Nn} \std{V}_d$ we will describe (canonically) a stably
dominated type $p_c$.   We will show that any stably dominated type on $V$ equals $p_c$ for a unique $c \in \std{V}$.  
 It follows that $c \in \dcl(A)$ iff $p_c$ is $A$-definable.  We define $\std{V} = \liminv_{d \in \Nn} \std{V}_d$.

In this language, \propref{b3} states that  there exists
a pro-definable partial map $f: \G^n \to \std{V}$ (over $M$) and $\g \in \G^n$ such that with $c=f(\g)$, we have $\g \in M(a)$
and $a \models p_c | M(\g)$.  

Thus $tp(a/M)$ can be understood in terms of (i) $tp(\g/M)$ and (ii) an $M$-definable function $\G^n \to \std{V}$. \>{discussion}

\<{ex} \lbl{compose}  Let $r$ be an $A$-definable type, and let $f$ be an $A$-pro-definable function into $\std{V}$, with $\dom(f) \in r|A$.
For any $B$ with $A \leq B$, let $a \models r |B$, $p=p_{f(a)}$, $c \models p|B(a)$.  Show that
$p_B = tp(c/B)$ does not depend on the choices, and that there exists a unique $A$-definable type $p$ with $p|B = p_B$.
We will refer to this type as $\int_r f$.  \>{ex}

In particular, \propref{b3} and the discussion below it yield:

\<{ex} Any $M$-definable type on $V$ has the form $\int_r f$ for some $M$-definable type $r$ on $\G^n$, and
some $M$-definable partial map $f: \G^n \to \std{V}$.  \>{ex}

We will later improve this to   decomposition theorem over other bases:
 Every definable type on $V$ can be decomposed into a definable type over $\G^n$, and a germ of a definable function into $\std{V}$.

 \<{ex}  \lbl{b6}Let $M$ be a maximally complete model, and $\g  \in \G^n$.  Then $M(\g) = \dcl(M \union \{\g\})$ is
 algebraically closed.
\>{ex}

Hint:  Let $N$ be a model containing $M(\g)$, and with $\G(N)$ generated by $\g$ over $\G(M)$.   For any $a \in N$,
by \propref{b3}, $tp(a/M(\g))$ extends to an $M(\g)$-definable type.  In general if $e \in \acl(B)$ and $tp(e/B)$
extends to a $B$-definable type, show that $e \in \dcl(B)$.

\<{remark} \label{3.10}

 Even over a base $A$ consisting of imaginaries, if $p$ is a stably dominated $A$-definable type and, then it is dominated by some  $A$-definable function $f$ into a finite-dimensional $k$-vector space.  This follows from a general descent principle for stably dominated types
and the elimination of imaginaries we will prove later.

\>{remark}

\ssec{Definable modules}

We consider definable $K$-vector spaces $V \cong K^n$.   When working over a base $A$ we  will always assume $V$ has a basis of $A$-definable points; this can be taken as the definition, but in fact is automatic, at least over nontrivially valued subfields, by the following version of
Hilbert 90:

\<{lem}   Let $F$ be a nontrivially valued field.    If $V$ is an $F$-definable $K$-space then $V$ has a basis of $F$-definable points. \>{lem}

\prf  We may assume $F=\dcl(F) \meet K$.
In this case, $F^{alg}$ is a model, so $V$ has a basis of points of $V(F^{alg})$.  This basis lies in   $V(F')$ for some finite Galois
extension $F'$ of $F$.  Now the automorphism group of $F'/F$ in the sense of ACVF and of ACF coincide,
by \lemref{galois}. Hence the usual Hilbert 90 applies.  \eprf

\<{lem} \lbl{galois}
Let $T$ be any expansion of the theory of fields, $F$ a subfield of a model $M$ of $T$ with $\dcl(F)=F$.  Let $F' \leq M$ be a finite normal
extension of $F$.  Then every field-theoretic automorphism of $F'/F$ is elementary.  \>{lem}
\prf   Let $G$ be the set of automorphisms of $F'/F$ that are elementary, i.e. preserve all formulas.
Then $Fix(G) = dcl(F) = F$.  By Galois theory, $G=Aut(F'/F)$ in the field theoretic sense. \eprf

\bigskip

Let $Mod_V$ be the set of definable $\Oo$-submodules of $V$.     $\L \in Mod_V$ is {\em $g$-closed} if $\L$ intersects any 1-dimensional $K$-subspace $U \leq V$ in a submodule of the form $\Oo c$ or $U$ or $(0)$.
$\L$ is a {\em semi-lattice} if it is $g$-closed and generates $V$ as a $K$-space.  $\L$ is a {\em lattice} if it is
$M$-isomorphic to $\Oo^{\dim V}$.

Let $V^*$ be the dual space to $V$; we identify $V^{**}$ with $V$, and write $(u,v)$ for the pairing $V \times V^* \to K$.    For $\L \in Mod_V$,
let $\L^*  = \{v \in V^*: (\forall a \in \L)  (a,v) \in \Mm  \}$.   In class we considered a different notion, namely
$\L^*_c =  \{v \in V^*: (\forall a \in \L)  (a,v) \in \Oo  \}$.

\<{ex}  Let $\dim(V)=1$, and $\L \in Mod_V$.  Then $\L = \Oo c$ or $\L=V$ or $\L=(0)$ or $\L = \Mm c$ for some $c \in M$.
\>{ex}

\<{ex} \lbl{mod3}  \<{enumerate}
\item  $^*$ and $^*_c$ are weakly inclusion-reversing maps $Mod_V \to Mod(V^*)$.   We have
$\L^{**} = \L$, and if $\L$ is closed also $(\L^*_c)^*_c=\L$.
\item    Let $M \models ACVF$.  If $\L \in Mod_V(M)$ then $\L$ is $M$-isomorphic to
$K^l \times \Oo^m \times \Mm^n$ for some $l,m,n$.
\item  If $\L$ is closed, then $\L \cong K^l \times \Oo^m$ for some $l,m$.
\item  $\L^{*}_c$ is always closed.
\item  Define $\L_c = (\L^*_c)^*_c$.  Then $\L$ is the smallest closed $\Oo$-module containing $\L$.
$\L$ contains $\Mm \L$.
\>{enumerate}
\>{ex}

It follows from  \exref{mod3} (3)   that the elements of $Mod_V$ are uniformly definable.

\<{ex} \lbl{mod4}  Let $A$ be a valued field and let $e_1,\ldots,e_n$  be (imaginary) codes for modules,
$\G(A(e_,\ldots,e_n) = \dcl( A \union \{e_1,\ldots,e_n \}) \meet \G$.
  Then there exists a maximally complete
field $N$ with $A \leq N$, $e_1,  \ldots  \in \dcl(A)$ and $\G(N)=\G(A(e_,\ldots,e_n))$.

\rm Hint:   This reduces to the case $n=1$, so $e=e_1$ codes a submodule $\Lam$ of $K^n$.  We may assume $\Lam$
generates $K^n$, and the dual module generates the dual space; so $\Lam$ contains no nonzero subspace of $K^n$.
Let $\Lam_c$ be the smallest lattice containing $\Lam$.  By adding to $A$ a generic basis for $\Lam_c$, we may
assume $\Lam_c= \Oo^n$.   By \exref{mod3} (5), $\Mm^n \subseteq \Lam$.  So to define $\Lam$ over $A$
it suffices to define $\Lam / \Mm^n$, a subspace of $\k^n$.  This can be done with parameters from $\k$.
If $\a \in K$, and $a$ is a generic element of $\res \inv (\a)$, show that $\G(A(a)) = \G(A)$.
\>{ex}

 \<{ex} \lbl{mod5}  Let $\L$ be a semi-lattice in $V$.  For $a \in V$, show  that $\{-\val(c): ca \in \Lam \}$ has a unique maximal element $v_{\Lam}(a) \in \G$,
 unless $a \subseteq \Lam$; in the latter case write $v_{\Lam}(a)= \infty$.    Show that  $v=v_{\Lam}$ satisfies
 $v(a+b) \geq \min v(a),v(b)$ and $v(cb) = v(c)+v(b)$ for $a,b \in V, c \in K$.  If ${\Lam}$ is a lattice, then $(V,v_{{\Lam}})$ is a valued
vector space.  Conversely, given $v$ with the above properties, ${\Lam}_v = \{a: v(a) \geq 0 \}$ is a semi-lattice,
and alattice of $v(V \m (0)) \subseteq \G$.  \>{ex}

\ssec{Pro-definable structure on $\std{V}$ }

Let $V$ be an affine variety, $V \subseteq \Aa^n$.

Let $H_d$ be the space of polynomials in $n$ variables of total degree $\leq d$.  Let $LH_d$ be
set of semi-lattices in $H_D$.

Let $\std{V}$ denote the stably dominated types on $V$.  We define $J_d: \std{V} \to LH_d$ by
$$J_d (p) = \{f \in H_d:  (d_p v) ( f(v) \in \Oo) \}$$
$$J = (J_1,J_2, \ldots): \std{V} \to \Pi_d LH_d $$

\<{prop} \<{enumerate}
\item  $J$ is 1-1.
\item  The image of $J$ is a pro-definable set.
\item  In fact, the image of $J_d$ is a definable set.
\item  Let $f(v,u)$ be a polynomial in variables $(v,u)=(v_1,\ldots,v_n,u_1,\ldots,u_m)$, of $v$-degree $\leq d$.
There exists a definable function $h: LH_d \times \Aa^m \to \G$ such that for any $p \in \std{V}$
and $b \in \Aa^m$, if $f_b(v) = f(v,b)$ then $ \val f_b(v) = h( J_d(p),b) $ is in $p | \Uu$.  In other words,
 $\val f_b(v) $ takes a constant value on generic realizations of $p$,   this value as a function of $p$ factors through
 $J_d(p)$, and it is uniformly definable over $LH_d$.
\>{enumerate}
\>{prop}

\prf  It suffices to prove this for $V=\Aa^n$.    Let $\bar{\Lam} = (\Lam_d)_{d \in \Nn} \in \Pi_d LH_d$.

Define $P(\bar{\Lam}) =   \{ \val( f(x) ) = v_{\Lam_d} (x) \}: d \in \Nn, f \in H_d \}$, where $v_{\Lam_d}$ is as in  Exercise  \ref{mod5}.

Now check that $D=\{\bar{\Lam}: P({\bar{\Lam}} \hbox{ is consistent } \}$ is a countable intersection of definable sets.
If $P(\bar{\Lam})$ is consistent, it generates a complete type over $\Uu$ (denoted the same way);   type
is always in $\std{V}$.   Thus $J(\std{V}) = D$;
this gives (2).  Since $P(J(p))$ generates $p$, we have (1).  With this definition of $J$, (4) is clear:  $h(\Lam,b)=v_\Lam(f_b)$.

(3) is Theorem 3.1.1 in \cite{HL}; see a more explicit proof in the Appendix.  
\eprf

\<{section}{$\G$-internal subsets of $\std{V}$}

\<{defn}  A definable set $D$ (possibly in imaginary sorts) is {$\G$-internal} if (possibly over additional parameters)
there exists a definable $Y \subset \G^n$ and a surjective    definable map $Y\to D$.  Equivalently, there exists
an injective definable map $D \to \G^n$.
\>{defn}

The equivalence in the definition uses  elimination of imaginaries for $\G$ (an easy result.)  In fact over
one parameter from $\G$, there even exist definable sets of representatives for any definable equivalence relation.
Let $f: Y \to D$ be surjective.  Let $W$ be a definable set of representatives for the relation $f(y)=f(y')$.
Then $g: D \to W$ defined by $f(g(d))=d$ is a definable injective map.

We can call  $D$   {\em almost $\G$-internal} if there
exists a finite-to-one definable map $D \to \G^n$.
In fact by \exref{b6},
almost $\G$-internal definable sets are $\G$-internal.  For sets of lattices this can also be seen by noting that  the proof of
\propref{7.5} goes through for almost $\G$-internal sets, and that the conclusion implies $\G$-internality.

If $D$ is $A$-definable,
it will turn out that the implicit parameters in the definition of $\G$-internality can be taken to be in $\acl(A)$.

\<{lem}  \lbl{7.1} Let $D$ be a $\G$-internal subset of $K^n$.  Then $D$ is finite.
\>{lem}

\prf  It suffices to show that every projection of $D$ to $K$ is finite; so we may assume $n=1$.  If $D$ is infinite,
it contains an infinite closed ball; over additional parameters there is therefore a definable surjective
map $D \to \k$.  However if $Y \subseteq \G^n$ there can be no surjective map $Y \to \k$, by the orthogonality of $\k,\G$.   This contradiction shows that $D$ is finite.
\eprf

\<{lem}  \lbl{7.2} Let $D$ be a $\G$-internal  set of closed balls of equal radius in $\Oo$, i.e. $D \subseteq K / c \Oo$.  Then $D$ is finite.  \>{lem}

\prf   Let $D' = \union D$.  If $D$ is infinite then $D'$ contains a closed ball $d \Oo + e$ with $\val(d)<\val(c)$.
Now $x \mapsto \res (d \inv (x-e))$ maps $d \Oo +e$ onto $\k$, and factors through $D$.  We obtain a contradiction as
in \lemref{7.2}.   \eprf

\<{lem}  \lbl{7.3}  Let $M$ be a model, $\g \in \G^N$.  Then any $M(\g)$-definable closed ball has a point in $M$.       \>{lem}

\prf   An $M(\g)$-definable closed ball $b$ lies in some     $\G$-internal set    $D$ of  closed balls.
By
\lemref{7.3}, we may take $D$ to be      linearly ordered by inclusion.  The intersection of all elements of $D$
is a ball $b'$ defined over $M$, closed or open, but nonempty; as $M$ is a model, we can choose a point of
$b'$ over $M$.  \eprf

\<{lem}  \lbl{7.3.0} Any $\G$-internal set  $D$ of balls is the union of
a finite number of definable subsets, each linearly ordered by inclusion.
      \>{lem}

\prf  Here we refer to Prop. 2.4.4 of \cite{hhmcrelle}.  \eprf

We call a lattice $\L$ {\em diagonal} for a basis $(b_1,\ldots,b_n)$ if there exist $c_1,\ldots,c_n \in K$ with
$\L = \sum \Oo c_i b_i$.  In other words, $\L = \oplus_i \Lambda \meet Kb_i$

\<{prop}  \lbl{7.5} Let $M$ be a model.  Let $\L$ be an $M(\g)$-definable lattice in $V=K^n$.  Then $\L$ has an $M$-definable diagonalizing basis.   Moreover if $e_1,\ldots,e_n$ is the standard basis, we can
choose a diagonalizing basis of the form $U e$, where
 strictly lower triangular  matrix over $M$ \footnote{'strict' here means: 1's on the diagonal.}

\>{prop}

\prf  The case $n=1$ is trivial.  Let $V_1$ be a one-dimensional subspace of $V$, $\bV = V/V_1$, $g: V \to \bV$
the canonical homomorphism.   Choose $b_1$ such that $V_1 \meet \Lambda = \Oo b_1$.  Let $\bar{\L} = g \L$.  By induction, there exists an $M$-definable  basis
$\bar{b_2},\ldots,\bar{b_n}$ diagonalizing  $\bar{\L}$; so $\bar{\L} = \sum c_i \Oo \bar{b_i} $ for some $c_i$, with
$c_i \Oo$ defined over $M(\g)$.   Now $g \inv  ( c_i \bar{b_i}) \neq \emptyset$, and
$g \inv( \bar{b_i})$ is a coset of $V_1$, so
$c_i \inv \Lambda \meet g \inv( \bar{b_i})$ is a closed ball in $V_1$.  By \lemref{7.3} it has an $M$-definable point $b_i$.   Any element of $\Lambda$ may be written as $v_1 + a_2 c_2 b_2 + \ldots + a_n c_n b_n$, with $v_1 \in V_1$, $a_i \in  \Oo$.
 So $v_1 \in V_1 \meet \Lambda$.  Thus  $\Lambda = \oplus_i \Lambda \Oo b_i$.   \eprf

  We may write $\Lambda = \oplus_{i=1}^n \Oo \gamma_i Ue_i   =  U \oplus_{i=1}^n \Oo \gamma_i  e_i = U S_\g \Oo^n$,
  where $S_\g$ is the diagonal matrix $(\gamma_1,\ldots,\gamma_n)$.   

Given $\Lambda$, the matrix $US_\g$ is determined up to multiplication on the right by an element of $B_n(\Oo)$;
and $S_\g$ is determined by $US_\g$;  the image  of $S_\g$  in $D_n  / D_n(\Oo) =  \G^n$   depends on $\Lambda$ alone,
and we denote it $G(\Lambda)$.  (This corresponds to the composed
homomorphism $\gamma: B_n \to \G^n$, composition of $B_n \to B_n/U_n = D_n$ with the natural map $D_n \to \G^n$.)

\<{cor} \lbl{7.6}  let $D$ be a $\G$-internal set of lattices.  Then there exist a finite partition $D = \union_{i=1}^r D_i$
and  bases $b^1,\ldots,b^r$ such that each $\Lambda
\in D_i$ is diagonal in $b^i$.   The bases $b^i$ are strictly upper triangular.    The function $G$ defined in the paragraph above
is injective on each $D_i$.
 \>{cor}

\prf   As the matrix $U$ in the conclusion of  \propref{7.5} is defined over $M$, while $\Lambda$ varies
over a definable set,    compactness assures the existence of finitely many matrices $U_1,\ldots,U_r$ over $M$,
such that each $\Lambda \in D$ has the form $U_i S_\g \Oo^n$ for some $i \leq r$ and for $\g=G(\Lambda)$.
 Let $D_i =  U_i D_n\Oo^n$.  

\eprf 

\>{section} 

\<{section}{Definable types in ACVF}

Let $M$ be a model.
We say that $tp(a/M)$ is definable if there exists a (necessarily unique) $M$-definable type $p$ with
$tp(a/M) = p|M$.

\<{lem} \lbl{8.0}  If $tp(a/M)$ is definable, and $c \in \acl(Ma)$, then $tp(ac/M)$ is definable. \>{lem}

\prf  Let $\phi(xy) \in tp(ac/M)$ be a formula such that $\phi(a,y)$ has $m$ solutions, with $m$ least
possible.  Then $p(x)|N \union \phi(x,y)$ generates a complete type over any elementary extension $N$.  By \lemref{dom1}, this is a definable type.    \eprf

\<{lem}\lbl{8.1}  Let $A$ be any subset of $\Uu^{eq}$, i.e. any set consisting possibly of imaginary elements.
Let $V \subseteq K^n$ be an $A$-definable set.  Then there exists a definable type on $V$, over $\Uu$, with finite orbit under $Aut(\Uu/A)$.  \>{lem}

This comes as close as possible to saying that $p$ is $A$-definable; one cannot do better since $V$ might be finite, or may have a finite but nontrivial definable quotient.

\prf  By induction on $n$.  If $n=1$, $V$ contains finitely many balls, each with some finite
union of sub-balls missing.  The generic type of one of these balls will do.  For $n>1$, let $\pi: K^n \to K^{n-1}$   be the projection, and let $p'$ be a definable type on $V'=\pi(V)$ with finite orbit.  Let
$M$ be a model containing $A$, and let
$a \models p' | M$.  Let $p''$ be a definable type on $\pi \inv(a)$ with finite orbit under $Aut(\Uu/A(a))$.  So $p''$ is $A(a')$ definable, with $a' \in acl(A(a))$.  Let $a'' \models p'' | M(a,a')$.  By \lemref{8.0},
$tp(aa'/M)$ is definable, and hence $tp(aa'a''/M)$ is definable, so $tp(aa''/M)$ is definable, i.e.
equals $p|M$ for some definable type $p$.  The number of conjugates of $p$ is at most
the number of conjugates of $a'/A(a)$.  \eprf

Let $r$ be an $A$-definable type on $\G^n$.  By a {\em pro-definable function on $r$ into $\std{V}$} we mean a pro-definable function $f$
represented by a sequence of definable functions $f_i$, such that $\dom(f_i) \in r|A$ for each $i$.  

Let  $f$ be an pro-definable function on $r$  into $\std{V}$ with $\dom(f) \in r|A$, whose $p$-germ is defined over $A$.
 Recall the definition of  $\int_r f$  (\exref{compose}).  It depends on $f$ only through the $p$-germ of $f$, so
 that $\int_r f$ is an $A$-definable type.

\<{thm}  \lbl{8.5} Let $p$ be an $A$-definable type on a variety $V$.  Then there exist a definable type $r$ on $\G^n$ 
and a definable $r$-germ $f$ of  pro-definable maps   into  $\std{V}$, with $p=\int_r f$.\>{thm}

\prf

Let $M$ be a maximally complete model, containing $A$.
    
Let $c \models p |M$.  Let $\g=g'(c)$ be a basis
for $\G(M(c))$ over $\G(M)$; let $r'=g'_*p$.  
 
Now $tp(c/M(\g))$ is stably dominated, so it equals $q | M(\g)$
for some $q \in \std{V}$; we can write $q= f'(\g)$, with $f'$ an $M$-definable function into $\std{V}$.   
By definition, $p=\int_{r'} f'$.

\eprf
 The proof showed that  $r=g_* p$, where $g=\alpha^c \circ g'$.   In particular, the $r$-germ of
$g \circ f$ is the $r$-germ of the identity, i.e. $f$ is genericallly injective.   (We could also arrange this a posteriori.)

How canonical is the pair $(r,f)$?   

\begin{defn}
Consider pairs $(r,h)$ with $r$ a definable type  and $h$ a   definable function.  We say two such pairs
$(r,h) , (r',h')$ are {\em equivalent   up to generic reparameterization},  $(r,h) \sim (r',h')$,  iff there
exist definable functions $\phi,\phi'$ such that $\phi_* r = \phi'_* r'$, and for some definable $h''$,
$h= h'' \circ \phi$ and $h' = h'' \circ \phi'$.  

When $h'$ is generically injective, this is equivalent to the existence of a 
  a definable $\phi$ such that $r'=\phi_*r$ and 
  $h=h' \circ \phi $ as an $r$-germ.  

If $h$ is pro-definable, with target $X=\liminv X_k$ and $\pi_k: X \to X_k$ the defining maps,
  we say 
$(r,h) \sim (r',h')$ if $(r,\pi_k \circ h) \sim (r',\pi_k \circ h')$ for each $k$.
\end{defn}

\begin{lem}  The pair  $(r,f)$ is determined by $p=\int_r f$, up to generic reparameterization.  \end{lem}

\prf  Suppose $p=\int_r f = \int_{r'} f'$, with $r,r',f',f'$ defined over some $N$.  Let $\g models r | N$, $c \models f(\g) | N(\g)$.
So $c \models p|N$.   Since also $p=\int_{r'}f'$, we may find $\g' \models r'|N$ such that $c \models f'(\g') | N(\g')$.  
By stable domination of
$p$, we have $\G(N(c)) \subset N(\g)$.  We claim that $\g \in \G(N(c)) $.   Let
    $\g''$ be a basis
for $\G(N(c))$ over $N$.  Then $ tp(c/N(\g''))$ extends to a stably dominated type $p''$ defined over $N(\g'')$.  By
orthogonality to $\G$ again, $p''$ implies a complete type over $N(\g)$, namely $tp(c/N(\g))=p$.  It follows that $p=p''$
is based on $N(\g'')$, and so by generic injectivity of $f$ we have $\g' \in \dcl(N(\g''))$.  Thus $N(\G(N(c)))= N(\g)$
and    similarly $N(\G(N(c))) = N(\g')$.  So $N(\g)=N(\g')$.   Moreover $h(\g),h'(\g')$ are stably dominated types 
based on $N(\g)$ and with the same restriction to this base, namely $tp(c/N(\g))$;  so $h(\g)=h'(\g')$.  
Let $\phi$ be an invertible $N$-definable function such that
$\g' = \phi(\g)$; then $r' = \phi_* r$ and  as $h'(\g') =  h(\phi \inv(\g'))$, $h' = h \circ \phi \inv$.   \eprf

We will study 
this notion in the ACVF setting in  the next section, but we indicate now how it will go.
   We will see in \lemref{ii1} that after a possible reparametrization, one can find an   $A$-definable
function $G$ on $\std{V}$ such that $G \circ f$ is the identity germ on $r$.   (Basically this is the $0$-definable function $G$ of \corref{7.6}; we need $A$ only in order to find an affine patch $V'$ of $V$ and identify $\std{V'}$ with a sequence of lattices.)
This implies that $r=G_*p$ is $A$-definable, and also rigidifies $f$ so that reparameterization is no longer possible,
and the $r$-germ of $f$ is uniquely determined.  {\em Hence with these choices we find an $A$-definable $r$ and 
a function $f$ with $A$-pro-definable germ.}  We can even use \lemref{cc7} to make $r$, if we wish, $0$-definable; this requires an additional
  reparamterization by a certain $A$-definable translation.
 
 \begin{remark}   
   Though the $r$-germ of $f$ can be chosen to be $A$-pro-definable, it is not always possible to find an $A$-(pro)definable 
$f$.  For instance for the generic type of an $A$-definable open ball without an $A$-definable sub-ball, this is the case.
This phenomenon is responsible for much of the subtlety in  the stability-theoretic study of ACVF.

The   function $G$ described above, inverting the germ $f$ on the left, cannot in general be take of the form $p \mapsto g_*p$ for any $A$-definable $g$. 
\end{remark}

\>{section}  

\<{section}{Imaginaries in ACVF}

Recall $B_n$ denote the group of invertible upper triangular matrices.   $U_n$ is the group of matrices in $B_n$
with $1$'s on the diagonal.  $D_n$ is the group of diagonal matrices, so that $B_n = D_n U_n$.

If $G$ is any algebraic subgroup of the group $GL_n$ of invertible $n \times n$ - matrices, $G(\Oo)$ denotes
the elements $M \in G$ such that $M, M \inv$ have entries in $\Oo$.

Let $S_n$ be the coset space $B_n / B_n(\Oo)$.  We will see below that any lattice in $K^n$ has a triangular basis.
Hence $B_n$ acts transitively on the set of lattices; and $B_n(\Oo)$ is the stabilizer of the standard lattice $\Oo^n$.
It follows that $B_n / B_n(\Oo)$ can be identified with the set of lattices in $K^n$.  (By a similar argument, so can $GL_n(K)/GL_n(\Oo)$.)

Let $\widetilde{GL_n}(\Oo)$ be the pullback of the stabilizer of a vector, under  the natural homomorphism $GL_n(\Oo) \to GL_n(\k)$.    Let $T_n$ be the coset space $GL_n / \widetilde{GL_n}(\Oo)$  We have a natural map $T_n \to S_n$.  
Given $b \in S_n$, viewed as a lattice $\Lambda$,
naming an element of $T_n$ is equivalent to choosing a point of $\Lambda/ \Mm \Lambda$.   
  Let $GG$ consist of the valued field sort $K$, along with the  sorts $S_n,T_n$.

Certain related imaginary sorts can be directly shown to be coded in the sorts $S_n,T_n$.

\<{lem} \lbl{direct} \<{enumerate}
\item  Any definable $\Oo$-submodule of $K^n$, as well  any coset of such a submodule of $K^n$, can be coded in $GG$.
\item  Any finite subset of $S_n \union T_n \union K^m$ is coded in $GG$.
\item Let $H$ be a  subgroup of $U_n$ defined by a conjunction
\[    H=\{a \in U_n:  \bigwedge_{i \leq j \leq n} \val(a_{ij}) \diamond_{ij} \alpha_{ij}  \}   \]
where $\alpha_{ij} \in \G_\infty$ and $\diamond$ denotes $\geq$ or $>$.     Let $A$ be a base structure
containing $\alpha_{ij}, i,j \leq n$.  Then any coset of $H$ is coded in $GG_A$ (i.e. for any coset $C$
of $H$ there exists $g \in GG^m$ such that $g$ is a canonical code for $C$ over $A$.)
\>{enumerate}
\>{lem}

 \prf We will not repeat the proofs of (1,2) from \cite{hhmcrelle}; (1) is rather straightforward, see 2.6.6;  (2) is Prop. 3.4.1 there.
  .

For (3), let $A_n$ be the $\Oo$-algebra of strict \footnote{'strict' here means:  0's on the diagonal} upper triangular matrices.   Let $J$ be the subalgebra defined by:
$\bigwedge_{i \leq j \leq n} \val(a_{ij}) \diamond_{ij} \alpha_{ij} $.  Then
   $H=1+J$.   We have $aH = bH $ iff $a=b(1+j)$ for some  $j\in J$ iff $aJ=bJ =: J'$ and $a+J'=b+J'$.  As $J'$ is an $\Oo$-module and
   $a+J'$ a coset,   (3) follows from (1).
 \eprf

\<{lem}  \lbl{ii1} Let $r$ be a definable type on a definable $D \subset \G^n$, $V=K^N$, and $h: D \to LV$
be an injective definable map.  Then $(r,h)/\sim$ has a canonical base in GG. \>{lem}

\prf   Let $U(\Lambda)$ be the maximal $K$-subspace contained in $\Lambda \in L$.
Say $\dim U(h(t))= d$.   $U(h(t))$ can be viewed as an element of a Grassmanian variety $Gr_{d} (V)$.  By \lemref{7.1},
 the image of $U(h(t))$ is finite.  Since $p$ is complete, the image is a single element $U$, i.e.  $U(h(t))=U_d$ for all $t \models r$.  Now $U$ is clearly an invariant
of $(r,h)/ \sim$.  We may work over a base where all $U$ are defined, and view
$h$ as a function $r \to L(V/U)$.   We may thus assume
$h(t)$ is a lattice for $t \models r$.

By \corref{7.6} there exists a triangular basis $b$ for $V$ such that $h(t)$ is diagonal in $b=(b_1,\ldots,b_n)$, for $t \models r$.  So $h(t) = \sum \Oo \g_i(t) b_i$ for certain definable
functions $\g_i : r \to \G$.  Let $\g=(\g_1,\ldots,\g_n)$.
(Recall  $\Oo \g $ denotes $\{x: \val(x) \geq \g \}$. )  
  We can replace $r$
by $\g_* r$ and $h$ by the function $h_b((s_1,\ldots,s_n)) = \sum \Oo s_i b_i$,
without changing the $\sim$-class.  So from now on we will consider only $h=h_b$ of this form.  Thus we need to code pairs $(r,b)$ up to $\sim$, where $(r,b) \sim (r',b')$
iff $(r,h_b) \sim (r',h_{b'})$.  Note that $h_b$ is injective on $\G^n$.

By \lemref{cc7} there exists $c \in \G^n$ such that $\a^c_* (r)$ is 0-definable; where $\a^c$ is translation by $c$.
Say $c=(c_1,\ldots,c_n), c_i =  -\val(e_i)$.  Let $e$ be the diagonal matrix $(e_1,\ldots,e_n)$.  Then
 $(r,b) \sim (\a^c_*(r),eb)$.  Replacing  $(r,b)$ by $ (\a^c_*(r),eb)$, we may assume $r$ is 0-definable.

Since $r$ is 0-definable,    $(r,b) /\sim$ is equi-definable with $b/\sim$, so
 we will now fix $r$ and consider the equivalence relation: $b \sim b'$ iff $(r,b) \sim (r,b')$

We view $b=(b_1,\ldots,b_n)$ as a matrix, with $b_i$ the $i$'th column.
Note   that $b,b'$ generate the same $\Oo$-lattice iff $b'_i = \sum c_{ij} b_j$ for some $c_{ij} \in \Oo$
and conversely, iff $b' = b N$ for some $N \in GL_n(\Oo)$.  Also, we have $(t_1b_1,\ldots,t_nb_n) ,bD_t$
generate the same $\Oo$-module,
where $D_t$ denotes any triangular matrix $(e_1,\ldots,e_n)$
with $\val(e_i)=t_i$.

Suppose $b \sim b'$.  So $r'=\phi_* r$ for some definable function
$\phi$ and $h_{b} = h_{b'} \circ \phi$.  Let $N=N(b,b')$ be the change of basis matrix, $bN(b,b') =b'$.
Then $N(b,b')$ is upper triangular.  Write $s=\phi(t)$.  Then $bD_t$ and 
 $b'D_s = b N D_s$  generate the same $\Oo$-module,  so $b D_t = bND_s  N'$  for some $N' \in GL_n(\Oo)$, or $D_s \inv N D_t \in   GL_n(\Oo)$.
Equivalently $D_s \inv N D_t \in   B_n(\Oo)$.
But $N$ is upper triangular,   $N(b,b') =  D(b,b') U(b,b')$, with $D(b,b')$ diagonal
and $U(b,b')$ strictly upper triangular.  It follows that $D(b,b')  =D_{t-s} \mod D(\Oo)$.  This holds
for $t \models r$; so $t-s$ is constant, i.e. $\phi(t)=t+c_0$, where $c_0 \in S(r) = \{c \in \G^n: \a^c_* r = r \}$.
Note that $S(r)$ is a definable subgroup of $\G^n$ (of the form  $E ( \G^l \times (0))$ for
some $l \leq n$ and some matrix $E$ with $\Qq$-coefficients.)  Let $S'(r)$ be the pullback of $S(r)$ to the group $D_n(K)$ of diagonal matrices.
Then $D(b,b') \in S'(r)$.  Moreover, since $D_s \inv N D_t \in   B_n(\Oo)$,   we have $U(b,b') \in D_t B_n(\Oo) D_t \inv$,
or $U(b,b') \in D_t U_n(\Oo) D_t \inv$..
Conversely, the argument reverses to show that if $D(b,b') \in   S'(r)$ and $U(b,b') \in D_t U_n(\Oo) D_t \inv$ for
generic $t \models r$, then $b \sim b'$.   Let $D_\nu = \{g \in U_n: (d_r t)(g \in D_t U_n(\Oo) D_t \inv ) \}$.
It is easy to see that this is one of the groups in \lemref{direct} (3), and hence coded in GG.

\eprf

\<{thm}  \lbl{ii2}  In the sorts GG, ACVF admits elimination of imaginaries.\>{thm}

\prf  By  \lemref{8.1}, \lemref{dd3} and \lemref{direct} (2), 
  it suffices to show that any
definable type $q$ on $V=\Aa^n$ has a canonical base in the sorts $GG$.  Now
$q$ has the form $ \int_r h$ where $r$ is a definable type on $\G^m$ and
$h: \G^m \to \std{V}$ is a definable map.  $q$ is equi-definable with the pair $(r,h)$
  up to generic reparameterization.

We have $h=(h_d)$, $h_d: r \to LH_d$, where $H_d$ is the space of polynomials in $n$ variables of degree $\leq d$.   Define $\sim_d$ as $\sim$ above.
For large enough $d$, $h_d$ is injective on a definable neighborhood of $r$.  \footnote{alternatively, for any $d$, we can factor out the kernel of $h_d$
and work with the pushforward $r_d$ of $r$.}
Clearly if $\si$
fixes $q$ then it fixes $(r,h_d)/\sim_d$ for each $d$; conversely if
$\si$ fixes  $(r,h_d)/\sim_d$ for large enough $d$, then it fixes the $q$-definition of
any given formula, so it fixes $q$.  Thus it suffices to code $(r,h_d)/\sim_d$ for each $d$.  This
was proved in \lemref{ii1}.  \eprf

\>{section} 

\begin{section}{Appendix}

We give here an effective description of the image of $\std{V}$ in the space of semi-lattices.  This description came out of a conversation with  Bernd Sturmfels.

Let $F$ be a valued field.
We will use Robinson's quantifier-elimination theorem in a two-sorted version, i.e. some variables range over $K$ and some range over the residue field $k$.   This follows easily from the one-sorted version:   
if $\phi(x_1,\ldots,x_n)$ is quantifier-free formula on $\Oo^n$, which is invariant under translation
by $\Mm^n$, then the solution set of $\phi$ can be viewed as a subset of $k^n$; and it is easy to see that this subset is constructible
(a Boolean combination of varieties.)   Note that  if $\phi(x_1,\ldots,x_n,y_1,\ldots,y_m)$ is $\Mm^n$-translation invariant in the $x$-variables,
with the $y$-variables fixed, then so is $(\exists y_1) (\exists y_2) \phi$, or any other sequence of quantifiers over the $y$-variables. 

\ssec{}   Let us take an affine variety 
$V = Spec(k[X_1,\ldots,X_m] / I)$.  Let $\std{V}$ be the stable completion.
 Let $H_d$ be the vector space of   polynomials of degree $\leq d$,  
$I_d = H_d \meet I$, $U_d=H_d / I_d$, and let $S(U_d)$ be the space of  semi-lattices in $U_d$.    

There is a natural map $r_{d,n}: \std{V} \to S(U_d)$.   Namely if $p$ is viewed as a semi-norm, $r_{d,n}(p)= \{f+I_d: p(f) \geq 0 \}$.

Given $\Lambda \in S(U_d)$, let $T=T(\Lambda)$ be the maximal $K$-space contained in $\Lambda$,
and let  $(f_1,\ldots,f_n)$ be an $\Oo$-basis for $\Lambda/T$.   
Let 
\[  R(\Lambda) = \{(\res f_1(c),\ldots,\res  f_n(c) ):  c \in V(K), v (f(c)) \geq 0 \hbox{ for } f \in \Lambda \} \] 

By Robinson's theorem, this is a constructible subset of $k^n$.   
 If we change the $\Oo$-basis, $R(\Lambda)$ changes by a linear transformation.

 \<{lem}  $\Lambda \in r_{d,n} (\std{V}) $ iff  $R(\Lambda)$ is not contained in a finite union of proper subspaces of $k^n$.    \>{lem}

\prf   First suppose $\Lambda \in r_{d,n} (\std{V}) $; say $\Lambda = r_{d,n}(p)$.  Suppose $R(\Lambda)$ is contained  in a finite union of proper subspaces of $k^n$; these subspaces and all data are defined over some model $M$.  Let $c \models p|M$;
let $f_1,\ldots,f_n$ be an $\Oo$-basis for $\Lambda$; then $ (\res f_1(c),\ldots,\res  f_n(c) ) \in R(\Lambda)$, 
so it must lie in one of the $M$-definable proper subspaces mentioned above; i.e. $\sum \alpha_i \res f_i(c) =0$, $\alpha_i \in k(M)$,
not all $0$.  Extend $(\alpha_1,\ldots,\alpha_n)$, viewed as an element $\chi_1$ of $(k^n)^*$,  to a basis $\bar{\chi_1},\ldots,\bar{\chi_n}$ of  $(k^n)^*$,  and lift to  a basis $\chi_1,\ldots, \chi_n$ of  dual lattice
$\Lambda^*$ of $\Lambda$.  Let $g_1,\ldots,g_n$ be the dual basis of $\Lambda$.  Then whenever $a \models p$, $g_1(a)$ has positive
valuation; say $\alpha = \val(c)$; it follows that $c \inv g_1 \in \Lambda$, but $c \inv \notin \Oo$, a contradiction.  


Conversely, assume   $R(\Lambda)$ is not contained in a finite union of proper subspaces of $k^n$.   Let $M$ be a maximally complete model over which $V,\Lambda$ are defined, let $f=(f_1,\ldots,f_n)$ be a basis for $\Lambda$ over $M$.  
Find  $c \in V$ be such that $f(c)=0$ for $f \in T$, $\val f_i(c) \geq 0$ for $i \leq n$ and $(\res f_1(c),\ldots,\res f_n(c)) $   does not lie in any proper $M$-definable subspace of $k^n$.  Let $\a$ be a basis
for $\G(M(c))$ over $\G(M)$;  so there exists a stably dominated type $p$ over $M(\a)$ with $c \models p |M(\a)$.  
It is clear that  $T \subseteq I(p) \meet H_d$ (where $I(p)$ is the kernel of the semi-valuation $p$) 
and $\Lambda  \subset r_{d,n}(p)$.   
We claim
that in fact, $r_{d,n}(p)=\Lambda$.   
  For suppose (e.g.) that $I_d =  I(p) \meet H_d$
but
 $r_{d,n}(p)$ is a bigger lattice $\Lambda'$.  The lattice $\Lambda'$ is defined over $M(\a)$ and so lies in a $\G$-parameterized
 family of lattices over $M$, so there exists a basis $g_1,\ldots,g_n$ of $H_d / T$ such that $\Lambda'$ is diagonal in this basis,
 i.e. $Q$ is generated by $c_1g_1,\ldots,c_ng_n$ for some $c_1,\ldots,c_n$.  
 The change-of-basis matrix $Q$ from $c_1g_1,\ldots,c_ng_n$ to $f_1,\ldots,f_n$ lies in $M_n(\Oo)$; if it is in $GL_n(\Oo)$,
 then the lattices are equal; if not, then 
some element $e$ of $\Lambda(M)$ lies in $\Mm \Lambda'$ but not in $\Mm \Lambda$.  As $e \notin \Mm \Lambda$,
we have $\res e(c) \neq 0$, otherwise $(\res f_1(c_m), \ldots,\res f_n(c_m))$ would lie in a proper subspace.  
 It follows that $\val e(c)=0$, and so   $p(e)=0$, contradicting $e \in \Mm \Lambda'$.
 
 \eprf
  
Now one can algorithmically decompose the constructible set $R(\Lambda)$ into irreducible, relatively closed sets and find their linear span; the condition of the lemma is that one of these spans should have dimension $n$.   This gives an effective description of the image of $\std{V}$.

\ssec{}   We have in general $\dim(R(\Lambda)) \leq \dim(V)$.     An important subset of the stable completion (denoted $V^\#$)  consists of the 
  strongly stably dominated points (see \cite{HL}).   In the present setting, a stably dominated type 
  on a variety $V$  is strongly stably dominated iff the residue field extension it induces has the same transcendence degree as the field extension it induces.   
  
Now if $\Lambda$ is a lattice with  $\dim(R(\Lambda))=\dim(V)$, then  $\Lambda$ is the image of  at most   a finite number $n(\Lambda)$ of elements $p$ of $\std{V}$, such that  for  $g_1,\ldots,g_n$ a basis of $\Lambda$, $M$ a model over which the data is defined,
and $c \models p|M$, $\res g_1(c),\ldots,\res g_n(c)$
are linearly independent over $\k(M)$.  These points $p$ all lie in  $V^\#$; 
and an upper bound on their number is easily given.   
This raises the hope of describing elements of $V^\#$ via a single tropical approximation.   But we have:  

\<{problem} Let $\Lambda$ be given, and assume $\dim(R(\Lambda))=\dim(V)$.   Determine $n(\Lambda)$ (or just whether $n(\Lambda)=1$) effectively.   
\>{problem}


\end{section}

 \end{document}